\documentclass[righttag,12pt]{amsart}
\usepackage{amssymb}
\usepackage{euscript}

\let\cal\mathcal

\newtheorem{prop}{Proposition}[section]
\newtheorem{theorem}[prop]{Theorem}
\newtheorem{cor}[prop]{Corollary}
\newtheorem{lem}[prop]{Lemma}

\newtheorem{example}[prop]{Example}

\theoremstyle{definition}

\newtheorem{ack}{Acknowledgments} 
\theoremstyle{remark}
\newtheorem{rem}[prop]{Remark} 

\newcommand{\bbN}{{\mathbb{N}}}

\newcommand{\bbR}{{\mathbb{R}}}

\newcommand{\calA}{{\cal{A}}}

\newcommand{\calE}{{\cal{E}}}

\newcommand{\al}{\alpha}

\newcommand{\de}{\delta}

\newcommand{\f}{\varphi}

\newcommand{\Sig}{\Sigma}

\renewcommand{\span}{\operatorname{span}}
\newcommand{\supp}{\operatorname{supp}}

\newcommand{\card}{\operatorname{card}}

\newcommand{\disp}{\displaystyle}

\newcommand{\lb}{\label}

\newcommand{\wtw}{if and only if }
\newcommand{\ceo}{conditional expectation operator}
\newcommand{\ceos}{conditional expectation operators}
\newcommand{\wceo}{weighted conditional expectation operator}
\newcommand{\wceos}{weighted conditional expectation operators}

\newcommand{\DEF}{\buildrel {\mbox{\small def}}\over =}

\begin{document}

\title{A disjointness type property 
\\ of conditional expectation operators}

\author{Beata Randrianantoanina$^*$}\thanks{$^*$Participant, NSF Workshop
in Linear Analysis and Probability, Texas A\&M University}

\address{Department of Mathematics and Statistics
\\ Miami University \\Oxford, OH 45056}

 \email{randrib@muohio.edu}



\begin{abstract} We give a characterization of conditional expectation
operators through a disjointness type property similar to band preserving 
operators. We say that the operator $T:X\to X$ on a Banach 
lattice $X$ is semi band preserving if and only if for all $f, g 
\in X$, $f \perp Tg$ implies that $Tf \perp Tg$. We prove that when
 $X$ is a purely atomic Banach lattice, then an operator $T$ on $X$ 
is a weighted conditional expectation operator
 if and only if $T$ is  semi band preserving.
\end{abstract}

 \maketitle

\section{Introduction}

In this note we study two abstract disjointness type conditions 
which are satisfied by all conditional expectation operators on 
Banach lattices. There is an extensive literature devoted to 
finding conditions which characterize conditional expectation 
operators and an extensive literature studying disjointness preserving
and band preserving operators.  However, as far as we know, to date there 
have  been no attempts to characterize conditional expectation 
operators through a property related to disjointness.

Of course, conditional expectation operators are never 
disjointness preserving yet alone band preserving.  However they do 
preserve some bands, namely they satisfy the following 
disjointness type condition:
\begin{equation}\lb{D}\tag{SBP}
f \perp Tg  \Longrightarrow Tf \perp Tg   \ \ \forall f,g \in X,
\end{equation}
(here $X$ is a Banach lattice and $T$ is a linear operator on 
$X$). 

Note that the condition \eqref{D} is   a weakening of the condition which
defines band preserving operators. Recall that a linear operator $T$ on
a Banach lattice $X$ is called {\it band preserving} if $TB\subset B$
for every band $B\subset X$. Thus $T$ is band preserving \wtw
one of the two following equivalent conditions is satisfied.
\begin{equation}\tag{BP1}
f \perp g  \Longrightarrow Tf \perp g  \ \ \forall f,g \in X,
\end{equation}
\begin{equation}\tag{BP2}
f \vartriangleleft g \Longrightarrow  Tf   \vartriangleleft g   \ 
\ \forall f,g \in X.
\end{equation}
(We use notation $f \vartriangleleft g$ to mean that $f$ belongs to 
 a band generated by $\{g\}$.)

Thus condition \eqref{D}  is the same as  $(BP1)$ with the 
additional constraint that $g$ belongs to the range of $T$. Hence,
clearly $(BP1)$ implies \eqref{D} and $(BP1)$ and \eqref{D} are
equivalent if $T$ is surjective. Conditional expectation operators are
our principal examples of non-band preserving operators which do
satisfy \eqref{D}.

We will say that an operator $T$ is {\it semi band preserving} if 
$T$ satisfies \eqref{D}. 
Our main result (Theorem~\ref{LT} and Corollary~\ref{charceo}) asserts
that when
 $X$ is a purely atomic Banach lattice, then an operator $T$ on $X$ 
is a \wceo ~\wtw $T$ is  semi band preserving.

Further, we study a condition  which arises from the weakening of 
$(BP2)$ by adding the 
 constraint that $g$ belongs to the range of $T$, 
similarly as in the definition of semi band preserving operators.
Namely we consider 
\begin{equation}\tag{SCP}\lb{C}
f \vartriangleleft Tg \Longrightarrow  Tf   \vartriangleleft Tg   \ 
\ \forall f,g \in X.
\end{equation}
We will say that an operator $T$ is {\it semi containment preserving} if 
$T$ satisfies \eqref{C}. 
It is clear that all surjective semi containment 
preserving operators are band preserving. It is also easy to see that 
all \ceos ~are  semi containment 
preserving   but not band preserving.
In contrast to the fact that $(BP1)$ and $(BP2)$ are equivalent,
conditions \eqref{D} and \eqref{C} are independent in general
(see Examples~\ref{ex1} and \ref{ex2}). However if Banach lattice $X$
is purely atomic then it follows from our characterization of semi
band preserving operators that all semi
band preserving operators are semi
containment preserving (see Corollary~\ref{eq}).
It is easy to construct on almost all Banach lattices a semi
containment preserving operator $T$ so that $T$ is not semi band preserving,
one can even find projections with this property (see Example~\ref{ex2}).
However we prove (Theorem~\ref{charscp} and Corollary~\ref{corscp}) that if $X$ is a strictly monotone purely atomic Banach
lattice and $P$ is a  projection of norm one on $X$ then $P$ is a  
\wceo ~\wtw $P$ is semi
containment preserving. (Thus, in particular, semi
containment preserving projections of norm one on 
 strictly monotone purely atomic Banach
lattices are   semi band preserving.)

We finish these general remarks about semi band preserving and semi
containment preserving operators by recalling a pair of conditions which
are very similar to \eqref{D} and \eqref{C}. Let $X$ denote a vector lattice
and $T$ be a linear operator on $X$. Consider:
\begin{equation}\tag{DP}
f \perp g  \Longrightarrow Tf \perp Tg  \ \ \forall f,g \in X,
\end{equation}
\begin{equation}\tag{$\beta$}
f \vartriangleleft g \Longrightarrow  Tf   \vartriangleleft Tg   \ 
\ \forall f,g \in X.
\end{equation}

Condition $(DP)$ is the well-known condition defining disjointness
preserving operators, and
  condition $(\beta)$ has been recently identified  by 
Abramovich and Kitover \cite{AK2000} as the condition equivalent to
the fact that $T^{-1}$ is disjointness
preserving (provided that $T$ is bijective and $X$ has sufficiently
many components).
Abramovich and Kitover \cite{AK2000} showed that in general 
 conditions $(DP)$ and 
$(\beta)$ are independent, but if $T$ is a continuous (or just regular)
linear operator between normed vector lattices then $(DP)$ implies 
$(\beta)$ and if $X$ is a Banach lattice lattice and $T$ is bijective then
$(DP)$ is equivalent to 
$(\beta)$.

\begin{ack} I wish to express my thanks to Professors Y. Abramovich
and A. Schep for their valuable remarks on preliminary versions of
this paper.
\end{ack}

\section{Preliminaries}

We use standard lattice and Banach space notations as may be found 
e.g. in \cite{LT1,LZ,MN}.
Below we recall basic definitions  that we use.

A {\it band} in a Banach lattice $X$ is a closed subspace $Y\subseteq X$
for which $y\in Y$ whenever $|y|\le |x|$ for some $x\in Y$ and so that
whenever a subset of $Y$ possesses a supremum in $X$, this supremum is a 
member of $Y$.
An element $u$ in a Banach lattice $X$ is called an {\it atom} 
if it follows from $0\neq v\le u$ that $v=u$. $X$ is called a
{\it purely atomic Banach lattice} if it is the band generated by its
atoms. Examples of purely atomic Banach lattices include $c_0,\ c,\ \ell_p\
(1\le p \le \infty)$ and Banach spaces with 1-unconditional bases.
A Banach lattice $X$ is called   {\it nonatomic} 
if it contains no atoms.

For an element $u$ in a Banach lattice $X$, an element $v\in X$ is
said to be a {\it component of $u$} if $|v|\wedge |u-v|=0$.
A lattice $X$ is called {\it essentially one-dimensional} if for any 
two non-disjoint elements $x_{1},x_{2} \in X$ there exist non-zero 
components $u_{1}$ of $x_{1}$ and $u_{2}$ of $x_{2}$ such that 
$u_{1}$ and $u_{2}$ are proportional. This class of lattices is 
strictly larger than purely atomic lattices and does include some 
nonatomic lattices, see \cite[Chapter~11]{AKmemoir}.

A Banach lattice $X$ is called {\it strictly monotone} if for all 
elements $x,  y$ in $X$ with $x, y >0$ we have $\|x+y\|>\|x\|$.

In this note we will mainly consider Banach lattices of 
(equivalence classes of) functions on  a $\sigma$-finite measure 
space ($\Omega, \Sigma,\mu$) which are subspaces of $L_1(\Omega, 
\Sigma, \mu) + L_{\infty}(\Omega, \Sigma, \mu)$.

By the Radon-Nikodym Theorem for each  $f \in L_1 (\Omega, 
\Sigma, \mu) + L_\infty (\Omega, \Sigma, \mu)$ and for every 
$\sigma$-subalgebra $\calA$ of $\Sigma$ so that $\mu$ restricted to 
$\calA$ is $\sigma$-finite (i.e. so that $\calA$ does not have atoms of 
infinite measure) there exists a unique, up to equality a.e., 
$\calA$-measurable locally integrable function $h$ so that
$$
\int_\Omega g h d \mu = \int_\Omega g f d \mu
$$
for every bounded, integrable and $\calA$-measurable  function $g$ on 
$\Omega$.  The function $h$ is called the {\it conditional 
expectation of $f$ with respect to $\calA$} and it is usually denoted 
by $\calE (f|\calA)$.  The operator $\calE(\cdot | \calA)$ is 
called the {\it \ceo  ~generated by} $\calA$.  Sometimes, particularly 
when $(\Omega, \Sigma, \mu)$ is purely atomic, $\calE ( \cdot 
| \calA)$ is also called an {\it averaging operator}.  When $X$ is a 
purely atomic Banach lattice with a basis $\{e_i\}_{i \in \bbN}$ 
then averaging operators on $X$ have the following form:

The $\sigma$-finite $\sigma$-subalgebra $\calA$ is generated by  a family
of mutually disjoint  finite subsets 
of $\Bbb{N}$, $\{A_j\}^\infty_{j=1}$, and for all  $x= 
\sum\limits^{\infty}_{i=1} x_i e_i$  the conditional 
expectation $\calE (x|\calA)$ is defined by:
$$
\calE (x | \calA) = \sum^{\infty}_{j=1} \left(\frac{1}{\card(A_j)}  
\sum_{n \in A_{j}} x_n\right) (\sum_{n \in A_j} e_{n}).
$$

Conditional expectation operators have been extensively  studied 
by many authors since 1930s, for one of the most recent 
presentations of the subject see \cite{Abook}.  One of the main 
directions in the research concerning conditional expectation 
operators is to identify a property or properties of an operator 
$T$ that guarantee that $T$ is a conditional expectation 
operator, see \cite{DHP90}.

Let $X$ be a Banach lattice of functions on $(\Omega, \Sigma, \mu)$
and let  $k \in L_1(\Omega,\Sigma, 
\mu) + L_{\infty}(\Omega, \Sigma, \mu)$, $ w\in X'$.  Then $\calE 
(wf|\calA)$ is well defined for all $f \in X$.
Assume in addition that $k \calE (w f | \calA)\in X$  for all 
$f \in X$ and put
$$
Tf = k \calE (wf|\calA).
$$
Thus defined operator $T$ is called a  {\it weighted conditional 
expectation operator}.  Note that when $X$ is a purely atomic 
Banach lattice or when $\calA$ is a $\sigma$-subalgebra of $\Sigma$ 
generated by a family of mutually disjoint sets 
$\{A_j\}^{\infty}_{j=1}$ of finite measure on $\Sigma$ then 
weighted conditional expectation operators on $X$ have the 
following form:
\begin{equation}\lb{ceo}
Tf = \sum^\infty_{j=1} \langle \psi_j,  f \rangle u_j
\end{equation}
where $\{\psi_j\}^\infty_{j=1} \subset X'$ and 
$\{u_j\}^\infty_{j=1} \subset X$ are so that for all $j$, $\supp 
 \psi_j \subset A_j$ and $\supp u_j \subset A_j$.

Recall that when  $X$ is a space of (equivalence classes of) 
functions on $(\Omega, \Sigma, \mu)$ then $\supp f$ is the minimal 
closed subset of $\Omega$ so that $f(t)=0$ for a.e. $t \in \Omega \setminus
\supp f$.

Note that a weighted  conditional expectation operator is a 
projection if and only if $\calE (kw|\calA)$ is the function 
constantly equal to $1$, in case when $\mu$ is a finite measure, 
or \wtw
$$
\langle\psi_j ,u_j\rangle=1 \ \text{for all} \  j
$$
in case when $\calA$ is a $\sigma$-subalgebra  of $\Sigma$ generated 
by a family of mutually disjoint sets $\{A_j\}_{j=1}^\infty$ (i.e. when
$T$ has form \eqref{ceo}).

\section{Definitions of semi band 
preserving and semi containment preserving operators}

Let $X$ is a Banach lattice and $T$ be a linear operator on 
$X$. As discussed in the Introduction we are interested in the following
two conditions:
\begin{equation}\tag{SBP}
f \perp Tg  \Longrightarrow Tf \perp Tg   \ \ \forall f,g \in X,
\end{equation}
\begin{equation}\tag{SCP}
f \vartriangleleft Tg \Longrightarrow  Tf   \vartriangleleft Tg   \ 
\ \forall f,g \in X.
\end{equation}
 We will say that an operator $T$ is {\it semi band preserving} if 
$T$ satisfies \eqref{D} and we will say that  $T$ is 
{\it semi containment preserving} if 
$T$ satisfies \eqref{C}.

It is easy to see that  all conditional expectation operators and \wceos
~are both semi band preserving and semi containment preserving.

Conditions \eqref{D}  and \eqref{C}  are weakenings of
conditions  $(BP1)$ and 
$(BP2)$ (respectively) which define band preserving operators, but
in contrast to the fact that conditions  $(BP1)$ and 
$(BP2)$ are always equivalent, in general
conditions \eqref{D}  and \eqref{C}  are independent of each other,
as the following
 two simple examples demonstrate.

\begin{example} \lb{ex1}
Let $X$ be a  Banach lattice of functions on $[0,1]$ such that 
the constant function $\f_1={\mathbf 1}=\chi_{[0,1]}$, and the 
function $\f_2$  defined by $\f_2(t)=t$ if $t \in 
[0,\frac{1}{2}],$ and $ \f_2(t)=0$ if $t \in ({1}/{2},1]$, belong 
to $X$ and there exist functionals
$\psi_1, \psi_2 \in X'$   with $\supp 
\psi_1\cup\supp\psi_2\subseteq[0,1/2]$. 
Then there exists a linear operator $T$ on $X$ which is 
semi band preserving but not semi containment preserving.
\end{example}

\begin{proof}[Construction] 
Define for all $f\in X$:
$$Tf = \langle\psi_1,f\rangle \f_1 + \langle\psi_2,f\rangle\f_2.$$

Then the operator $T$
 is semi band preserving. Indeed, $f\perp Tg$ implies that either
$f=0$ or $\supp Tg\subset [0,1/2]$ and $\supp f \subset [1/2,1]$. 
But then $Tf=0$ so $Tf\perp Tg$.

However $T$ is not semi containment preserving. Indeed, let $f, g 
\in X$ be such that $\langle\psi_1,f\rangle =0$,
$\langle\psi_1,g\rangle \ne 0$ and $\supp g\subset [0,1/2]$.
Then $Tf =   \langle\psi_2,f\rangle\f_2$ and so $\supp Tf= [0,1/2]$.
On the other hand, $\supp Tg= [0,1]$ since $\langle\psi_1,g\rangle \ne 0$.
Thus $g \vartriangleleft Tf$ but $ Tg   \not\vartriangleleft Tf$.
\end{proof}

\begin{example} \lb{ex2}
Let $X$ be any  Banach lattice which contains nonzero elements $f_1, f_2$
with $f_1\perp f_2$. Then there exists 
 a semi containment preserving operator $Q$ on X which is not 
semi band preserving. Moreover $Q$ can be chosen to be a projection
and if $X$ is not strictly monotone then $Q$ can be chosen to be a projection
of norm one.
\end{example}
\begin{proof}[Construction] 
Let  
$\psi $ be a functional  on $X$ so that $\langle\psi,f_1\rangle\ne 0$
and $\langle\psi,f_2\rangle\ne 0$. Define for all $f\in X$:
$$Qf = \langle\psi,f\rangle f_1.$$
Then  $Q$ is trivially semi containment preserving since the range of $Q$
is one dimensional. However
$Q$ is not 
semi band preserving since $f_2\perp Qf_1$, but $Qf_2\not\perp Qf_1$.

 Moreover if $\langle\psi,f_1\rangle=1$ then $Q$ is a projection.
Further if $X$ is not strictly monotone, then it is possible to chose 
$f_1\perp f_2$, $f_2\ne 0$,  so that $\|f_1 +f_2\|=\|f_1\|=1$ and  
 $\psi\in X'$ so that $\langle\psi,f_1\rangle=1$,
$\langle\psi,f_2\rangle\ne 0$ and $\|\psi\|=1$, which will result
in $Q$ being a projection of norm one.
\end{proof}

\section{Semi band preserving operators}

Our next goal is to characterize \wceos ~on purely atomic lattices as 
semi band preserving 
operators. 

In the following
 $X$ will be a Banach lattice of (equivalence classes of) real valued functions on a measure 
space $(\Omega,\Sigma,\mu)$. For any  linear 
operator $T : X \to X$ denote
$$
\Sigma_T = \{A \subset \Omega : \exists f \in X \ with \ \supp(Tf) 
= A\}.
$$
Here and in the following all set relations are considered modulo sets of 
measure zero.

We start with a simple lemma, which we formulate here for easy 
reference.

\begin{lem}\lb{sum} ~

\begin{itemize}
\item[(1)]If $A, B \in \Sigma_T$,  then $A\cup B \in \Sigma_T$.
\item[(2)]If $\{A_{j}\}_{j \in \bbN} \subset \Sigma_T$ is a family of
mutually disjoint sets, then $\bigcup ^{\infty}_{j = 1}
 A_{j} \in \Sigma_T$.
 \end{itemize}
\end{lem}

\begin{proof}
These facts are   immediate. For (1), let $f, g$ be concrete 
representations of functions in $X$ so 
that $\supp Tf=A$ and $ \supp Tg=B$. Define
\begin{equation*}
  h(t)\DEF
  \begin{cases}
    {\disp \frac{Tf(t)}{Tg(t)}} & \text{if}\ \  Tg(t)\ne 0, \\
    0 & \text{if}\ \  Tg(t)= 0,
  \end{cases}
\end{equation*}
and denote
$$V(h)=\{a\in\bbR : \mu(h^{-1}\{a\})>0\}.$$
Clearly $\card(V(h))\le \aleph_0$ and thus there exists 
$\al\in\bbR$, so that $-\al\notin V(h)$. It is easy to see that 
this implies that $\supp(T(f+\al g))= A\cup B$ (recall that all 
set relations are considered modulo sets of measure zero).

Part (2) is even quicker. Indeed let $\{f_{j}\}_{j \in \bbN}$ be
 a sequence of elements of $X$ such that $\|f_j\|=1$ and
 $\supp(Tf_{j}) = A_{j}$ for all $j \in \bbN$.
Then  $\sum ^{\infty}_{j=1} \ 2^{-j} f_{j}$ belongs to $X$ and, 
since sets $\{A_{j}\}_{j \in \bbN} $ are mutually disjoint
$$
\supp (T(\sum^{\infty}_{j=1}  2^{-j}f_{j})) = 
\bigcup^{\infty}_{j=1} A_{j},
$$
as desired.
\end{proof}

Denote $S_T\DEF \bigcup_{A\in\Sigma_T} A\subset \Omega$. Then, for 
each $f \in X$ we have
\begin{equation}
\supp (Tf) \subseteq S_T. \lb{supp}
\end{equation}

Now we immediately obtain:
\begin{prop}\lb{support P}
If $T$ is a semi band preserving  operator on $X$ then for every $f \in X$ with 
$\supp f \subseteq \Omega \setminus S_T$ we have $T f = 0$.
\end{prop}

\begin{proof}
Indeed, by \eqref{supp}, $\supp(Tf) \subseteq S_T$ so $f \perp Tf 
$. By \eqref{D}  we get $Tf \perp Tf$.  Thus $Tf = 0$.
\end{proof}

When the space $X$ is essentially one-dimensional we can deduce a further 
important property of semi band preserving operators. We have:

\begin{prop}\lb{complement}
Suppose that $X$ is essentially one-dimensional and $T$ is a semi 
band preserving operator on $X$.  If $A, B \in \Sigma_T$ and 
$A\subset B$, then $B\setminus A \in \Sigma_T$.
\end{prop}

\begin{proof}
Let $h, g \in X$ be such that $\supp Th = B$ and $\supp T g = A$. 
Since $A \subset B$ and $X$ is essentially one-dimensional there 
exists $C \subset A$ so that the components $(Th) \chi_{C}$ and 
$(Tg) \chi_{C}$ of $Th$ and $Tg$, respectively, are proportional. 
Let $\{A_{i}\}_{i \in I}$ denote the family of subsets of $A$ 
maximal with respect to the property that $(Th) \chi_{A_{i}}$ and 
$(Tg) \chi_{A_{i}}$ are proportional.  Then $\{A_{i}\}_{i \in I}$ 
are mutually disjoint and, by the essential one-dimensionality of 
$X$,
$$
A = \bigcup_{i \in I} A_{i}.
$$

Moreover, for each $i \in I$ there exists a scalar $a_{i} \ne 0$ 
so that
\begin{equation}\lb{C1}
(Th) \chi_{A_{i}} = a_{i} (Tg) \chi_{A_{i}}.
\end{equation}

Consider $g_{i} = h - a_{i}g$ for $i \in I$.  Then, by the 
maximality of $A_{i}$'s,
$$
\supp (Tg_{i}) = B \setminus A_{i}.
$$

Thus $B \setminus A_{i} \in \Sigma_T$ for all $i \in I$.  Moreover
$$
g \chi_{A_{i}} \perp Tg_{i}.
$$

Thus, by \eqref{D},
$$
T(g \chi_{A_{i}}) \perp Tg_{i}.
$$
That is, for all $i \in I$:
\begin{equation}\lb{C2}
\supp (T(g \chi_{A_{i}})) \subset A_{i}.
\end{equation}

But, since $\{A_{i}\}_{i \in I}$ are mutually disjoint
$$
Tg = \sum_{i \in I} \ T(g \chi_{A_{i}}),
$$
and, by \eqref{C2},
\begin{equation}\lb{C3}
(Tg) \chi_{A_{i}} = T(g \chi_{A_{i}}).
\end{equation}
Thus, by \eqref{C1} and \eqref{C3}, we get
\begin{eqnarray*}
(Th) \chi_{A} &=& \sum_{i \in I} \ (Th) \chi_{A_{i}} = \sum_{i 
\in I} \
a_{i}(Tg) \chi_{A_{i}}\\
&=& \sum_{i \in I} \ a_{i} T(g \chi_{A_{i}})=T(\sum_{i \in I} \ 
a_{i} g \chi_{A_{i}})=T(h \chi_{A}).
\end{eqnarray*}
Thus $(Th) \chi_{A} \in T(X)$. Hence
$$(Th) \chi_{B\setminus A} = Th- (Th) \chi_{A} \in T(X).$$
Thus $B \setminus A \in \Sigma_T$.
\end{proof}

\begin{rem} \lb{ess1}
Note that the above proof also shows that if $X$ is essentially 
one-dimensional and $T$ is a semi band preserving operator on $X$ 
  then the subspace $T(X)$ is essentially 
one-dimensional. We will prove a stronger result in Theorem~\ref{LT}.
\end{rem}

\begin{rem} \lb{LTex}
Proposition~\ref{complement} fails in general nonatomic Banach 
lattices. Indeed, let $T$ be the
  semi band preserving operator defined in Example~\ref{ex1}.
It is easy to see that $[0,1], [0,1/2]\in\Sig_T$ and 
$[1/2,1]=[0,1]\setminus [0,1/2]$ does not belong to $\Sig_T$.

Note that when $\psi_1$ and $\psi_2$ are positive then $T$ is 
positive, and when $\psi_i(\f_j)=\de_{ij} $ for $i,j = 1,2$, then 
$T$ is a projection. However it follows from 
\cite[Theorem~3.10]{DHP90} that when $T$ is an order continuous 
positive semi band preserving projection on a Banach lattice of 
functions on $[0,1]$   then $T$ satisfies the thesis of 
Proposition~\ref{complement}.
\end{rem}

By de Morgan Laws as a corollary of Lemma~\ref{sum} 
and Proposition~\ref{complement} 
we immediately obtain:
\begin{cor}\lb{intersection}
Suppose that $X$ is an essentially one-dimensional Banach lattice 
and $T$ is a semi band preserving operator on $X$. Then $T$ satisfies the
following two properties:
\begin{equation}\tag{$I1$}
A,B \in \Sigma_T\ \ \Longrightarrow \ \  A \cap B \in \Sigma_T;
\end{equation}
\begin{equation}\tag{$I2$}
\{A_{j}\}_{j \in \bbN} \subset \Sigma_T \text{ and }
A_1\supseteq
A_{2}\supseteq A_{3}\supseteq \dots\ \ \ \Longrightarrow \ \ \bigcap_{j \in \bbN} 
A_{j} \in \Sigma_T.
\end{equation}
\end{cor}

These properties allow us to give the full characterization
of semi band preserving operators on essentially 
one-dimensional Banach lattices. Namely we have:

\begin{theorem}\lb{LT}
Let $X$ be an essentially 
one-dimensional Banach lattice. Then an operator $T: X \to X$ is 
 semi band preserving  \wtw 
the range of 
$T$ is the linear
span of a collection of mutually disjoint elements 
$\{u_{j}\}_{j \in J}$ in 
$T(X)$ and $T$ is a \wceo, i.e. 
$T$ has the following form
for all $f$ in $X$:
\begin{equation}\lb{ceo1}
Tf=\sum_{j\in J} \langle \psi_j, f\rangle u_j,
\end{equation}
where  $\{\psi_j\}_{j\in J}$ are nonzero functionals on $X$ so that for
all $j\in J$ if $f\perp u_j$ then $\langle \psi_j, f\rangle=0$
 (see \eqref{ceo}).
\end{theorem}

\begin{proof}
It is not difficult to see that all \wceos ~are semi band preserving.

For the other direction,
let $\omega_{0} \in S_T \DEF \bigcup_{A\in\Sigma_T} A\subset 
\Omega$.  Then, by  $(I2)$ and Zorn's 
Lemma, among all $A \in \Sigma_T$ such that $\omega_{0} \in A$, 
there exists a set $A_{0}\in \Sigma_T$, minimal with respect to inclusion.

Next, we claim that the subspace of $T(X)$ consisting of those 
elements in $T(X)$ whose support is contained in $A_{0}$, is 
one-dimensional.  

Suppose for contradiction that there exist $f, g \in X$  such that 
$\supp Tf=A_{0}$, $\supp 
Tg \subseteq A_{0}$ and $Tf, Tg$ are linearly independent.  
Since $X$ is essentially one-dimensional there exist nonzero components
$u_1, u_2$ of $Tf, Tg $ respectively so that 
$$u_1=ku_2$$
for some scalar $k$. Clearly $\supp u_1 = \supp u_2$ and since 
$Tf, Tg$ are linearly independent
$$\supp u_1 = B\subsetneq A_{0}.$$
Consider $h=f-kg$. Then $Th= Tf-kTg$ and $C=\supp Th$ belongs to 
$\Sigma_T$ and 
$$\emptyset\ne C \subseteq A_{0}\setminus B \subsetneq A_{0}.$$
By Proposition~\ref{complement} we also have that $A_{0}\setminus C$
belongs to 
$\Sigma_T$.

Now $\omega_0$ belongs to one of the sets $C$ or  $A_{0}\setminus C$
 which 
contradicts the minimality of the set $A_{0}$.

It now follows immediately that there exist mutually disjoint 
elements $\{u_{j}\}_{j \in J}$ in $T(X)$ with minimal supports in 
$\Sigma_T$. Thus $T(X) = \overline{\span} \{u_{j}\}_{j \in J}$
and $T$ has the  form \eqref{ceo1} since $T$ is a linear operator.
Condition \eqref{D} implies that for all $j\in J$, if $f\perp u_j$,
since $u_j\in T(X)$, then $Tf\perp u_j$ and thus 
$\langle \psi_j, f\rangle=0$, as required in \eqref{ceo1}.
\end{proof}

\begin{rem}\lb{ex3}
The above proof  is very similar in spirit to the proof of the 
characterization
of the form of norm one projections in $\ell_p$, $1<p<\infty$,
\cite[Theorem~2.a.4]{LT1}.
\end{rem}
\begin{rem}
Theorem~\ref{LT} is not valid in general nonatomic lattices. 
The counterexample is very similar to Example~\ref{ex1}. 
Indeed let $X$ be any Banach lattice of 
functions on $[0,1]$ such that the constant function ${\mathbf 
1}$, and the function $\psi : [0,1] \to [0,1]$ defined by $\psi 
(t)=t$, belong to $X$.  Then $\span\{{\mathbf 1},\psi\} \subset 
X$ is 2-dimensional in $X$ and therefore it is complemented in 
$X$, i.e. there exists a projection $T : X \to X$ with $T(X)= \span 
\{\mathbf 1,\psi\}$.  But for every $g \in X$ we have $\supp Tg = [0,1]$.  
Thus $f\perp Tg$ implies $f=0$ and thus $T$ is 
trivially semi band preserving.
Clearly $T$ is 
not a \wceo. Further, note that every function in the range of $T$
has full support and hence $T$ is also trivially semi 
containment preserving.
\end{rem}

We finish this section with two immediate corollaries of Theorem~\ref{LT}.

\begin{cor}\lb{eq}
Let $X$ be an essentially 
one-dimensional Banach lattice. Then every semi band preserving operator
$T$ on $X$ is  semi containment preserving.
\end{cor} 

\begin{cor} \lb{charceo}
Let $X$ be a purely atomic Banach lattice. Then an operator $T$ on $X$ 
is a \wceo ~\wtw $T$ is  semi band preserving.
\end{cor} 

\section{Semi containment preserving projections}

In this section we obtain an analogue of our main result, Theorem~\ref{LT},
 for semi 
containment preserving operators. However, as Example~\ref{ex2} 
demonstrates,
on any  Banach lattice which contains nonzero elements $f_1, f_2$
with $f_1\perp f_2$  there exists 
 a semi containment preserving projection $Q$  which is not 
semi band preserving and thus is not a \wceo.
 Moreover if $X$ is not strictly monotone then such 
$Q$ can be chosen to be a projection of norm one.

Also an example described in Remark~\ref{ex3} demonstrates that in general 
nonatomic Banach lattices there may exist a semi containment 
preserving projection which is not a \wceo. Thus our characterization below
has natural restrictions. We prove:

\begin{theorem} \lb{charscp}
Let $X$ be an essentially 
one-dimensional strictly monotone Banach lattice and let $P:X 
\to X$ be a   projection of norm one. 
Then  $P$ is  semi containment preserving \wtw
the range of 
$P$ is the linear
span of a collection of mutually disjoint elements 
$\{u_{j}\}_{j \in J}$ in 
$P(X)$ and $P$ is a \wceo, i.e. 
$P$ has the following form
for all $f$ in $X$:
\begin{equation}\lb{ceo2}
Pf=\sum_{j\in J} \langle \psi_j, f\rangle u_j,
\end{equation}
where  $\{\psi_j\}_{j\in J}$ are nonzero functionals on $X$ so that for
all $j\in J$, $\supp \psi_j\subseteq\supp u_j$, 
 $\langle \psi_j, u_j\rangle=1=\|\psi_j\|=\|u_j\|$ and 
$\langle \psi_j, u_i\rangle=0$ if $i\ne j$
 (see \eqref{ceo}).
\end{theorem}

\begin{proof}
As before we note that all \wceos ~are semi containment preserving,
so we just need to prove one implication in Theorem~\ref{charscp}.

Our method of proof depends on the following lemma:

\begin{lem} \lb{propC}
Suppose that $X$ is a strictly monotone (not necessarily essentially 
one-di\-men\-sion\-al) Banach lattice and $P:X 
\to X$ is a  semi containment preserving projection of norm one.  
Let $\{A_{j}\}_{j \in \bbN} \subset \Sigma_P$ with $A_{1} 
\supseteq A_{2} \supseteq \dots$.  Then
$$
\bigcap_{j \in \bbN} A_{j} \in \Sigma_P.
$$
\end{lem}

Using this lemma the proof of Theorem~\ref{charscp} is the same as
the proof of Theorem~\ref{LT}. Indeed, Lemma~\ref{propC}
states that when $X$ and $P$ satisfy assumptions of Theorem~\ref{charscp}
then $P$ has property $(I2)$ from Corollary~\ref{intersection}.
Thus, following the proof of Theorem~\ref{LT} word for word, we get 
that there exist mutually disjoint 
elements $\{u_{j}\}_{j \in J}$ in $P(X)$ so that
 $P(X) = \overline{\span} \{u_{j}\}_{j \in J}$
and $P$ has the  form \eqref{ceo2} since $P$ is a linear operator.
Condition~\eqref{C} implies that for all $j\in J$, 
 $\supp \psi_j\subseteq\supp u_j$, and since $P$ is a projection of norm one
we have
 $\langle \psi_j, u_j\rangle=1=\|\psi_j\|=\|u_j\|$ and 
$\langle \psi_j, u_i\rangle=0$ if $i\ne j$, as required in \eqref{ceo2}.
\end{proof}

\begin{proof}[Proof of Lemma~\ref{propC}]
Since $\{A_{j}\}_{j \in \bbN} \subset \Sigma_P$, there exist 
$\{f_{j}\}_{j \in \bbN} \subset X$ so that $\supp Pf_{j}=A_{j}$.  
Denote $A=\bigcap_{j \in \bbN} A_{j}$ set $g=(Pf_{1})\cdot \chi_{A}$.  
Then $\supp g=A \subset \supp 
Pf_{j}$ for all $j \in \bbN$.  Thus, by \eqref{C},
$$
\supp Pg \subset \supp Pf_{j}
$$
for all $j \in \bbN$.  Hence
$$
\supp Pg \subset A.
$$
Denote $\supp Pg=B$.  Then
$$
(Pg) \cdot \chi_{A_{1}\setminus B}=0.
$$
Further
\begin{equation*}\begin{split}
Pf_{1} &=(Pf_{1})\cdot \chi_{A_{1}\setminus A}+(Pf_{1})\cdot 
\chi_{A}=(Pf_{1})\cdot \chi_{A_{1}\setminus A}+g,\cr Pf_{1} 
&=P(Pf_{1}) =P((Pf_{1})\cdot \chi_{A_{1}\setminus A})+Pg,\cr 
(Pf_{1})\cdot \chi_{A_{1}\setminus B} &=P((Pf_{1})\cdot 
\chi_{A_{1}\setminus A}) \cdot \chi_{A_{1}\setminus B} + (Pg)\cdot 
\chi_{A_{1}\setminus B}\cr &=P((Pf_{1})\cdot \chi_{A_{1}\setminus 
A})\cdot \chi_{A_{1}\setminus B}\ .
\end{split}\end{equation*}

Since $P$ has norm one we get:
\begin{equation*}
\begin{split}
\|(Pf_{1})\cdot \chi_{A_{1}\setminus B}\| &= \|P((Pf_{1})\cdot 
\chi_{A_{1}\setminus A})\cdot \chi_{A_{1}\setminus B}\|
\leq\|P((Pf_{1})\cdot\chi_{A_{1}\setminus A})\|\\
 &\leq\|(Pf_{1})\cdot
\chi_{A_{1}\setminus A}\|.
\end{split}\end{equation*}

Since $X$ is strictly monotone and $\supp Pf_{1}=A_{1}$ we 
conclude that
$$
A_{1}\setminus B \subseteq A_{1}\setminus A.
$$
Since $B \subset A$, we get that
$$
A=B=\supp Pg.
$$
Thus $A \in \Sigma_P$, as desired.
\end{proof}

We finish this section with an immediate corollary of 
Theorem~\ref{charscp} similar to Corollary~\ref{charceo}.

\begin{cor} \lb{corscp}
Let $X$ be a purely atomic Banach lattice and let $P:X 
\to X$ be a   projection of norm one.  Then an operator $P$ 
is a \wceo ~\wtw $T$ is  semi containment preserving.
\end{cor} 


\end{document}